\documentclass{article}
\usepackage{amsfonts}
\usepackage{enumerate}
\usepackage{amsmath,amssymb}
\usepackage[table]{xcolor}

\usepackage{float}

\newtheorem{thm}{Theorem}[section]
\newtheorem{lem}[thm]{Lemma}
\newtheorem{result}[thm]{Result}

\newtheorem{rmk}[thm]{ Remark}
\newtheorem{defi}[thm]{Definition}

\begin{document}

\vspace*{0.5cm}

\begin{center}
{\Large On automorphism groups of a biplane (121,16,2)}
\end{center}

\vspace*{0.5cm}

\begin{center}

Dean Crnkovi\'c \\
({\it\small E-mail: deanc@math.uniri.hr})\\[3pt]
Doris Dumi\v ci\'c Danilovi\'c\\
({\it\small E-mail: ddumicic@math.uniri.hr})\\[3pt]
Sanja Rukavina \\
({\it\small E-mail: sanjar@math.uniri.hr})\\[3pt]
 
\end{center}

\begin{center}
Department of Mathematics, University of Rijeka\\
Radmile Matej\v ci\'c 2, 51000 Rijeka, Croatia 
\par\end{center}

\bigskip

\vspace*{0.5cm}

\begin{abstract}
The existence of a biplane with parameters $(121,16,2)$ is an open problem.
Recently, it has been proved by Alavi, Daneshkhah and Praeger that the order of an automorphism group of a of possible biplane ${\mathcal D}$ of order $14$ divides $2^7\cdot3^2\cdot5\cdot7\cdot11\cdot13$.
In this paper we show that such a biplane do not have an automorphism of order $11$ or $13$, and thereby establish that $|Aut({\mathcal D})|$ divides $2^7\cdot3^2\cdot5\cdot7.$ 
Further, we study a possible action of an automorphism of order five or seven, and some small groups of order divisible by five or seven, on a biplane with parameters $(121,16,2)$.
\end{abstract}

\vspace*{0,5cm}



\section{Introduction}\label{intro}

We assume that the reader is familiar with the basic facts of design theory. 
For background reading we refer the reader to \cite{bjl} and \cite{crc}. 

One of the main problems in design theory is a classification of certain combinatorial designs with given parameters. 
An incidence structure ${\mathcal D} =( {\mathcal P},{\mathcal B},{\mathcal I})$, with point set ${\mathcal P}$, 
block set ${\mathcal B}$ and incidence ${\mathcal I}$ is called a $t$-$(v,k,\lambda)$ design, if $|{\mathcal P}|=v$, every block
$B \in {\mathcal B}$ is incident with precisely $k$ points, and every $t$ distinct points are together incident with precisely $\lambda$ blocks.
A design is called symmetric if it has the same number of points and blocks. 
The repetition number of each point of a symmetric design is $k$ and every two blocks are together incident with exactly $\lambda$ points.
A classification of designs has been considered in detail in the monograph \cite{k-o}. In the case of symmetric designs complete classification is done just for a few parameter triples (see \cite{mr}). 

A symmetric $2$-$(v,k,2)$ design is called a biplane. Recently, an overview of known biplanes and their automorphisms is given in \cite{biplanes}. 
An automorphism of a design ${\mathcal D}$ is a permutation on ${\mathcal P}$ which sends blocks to blocks.
The set of all automorphisms of ${\mathcal D}$ forms its full automorphism group which is denoted by $Aut({\mathcal D})$. 
Biplanes with parameters $(7,4,2)$, $(11,5,2)$, $(16,6,2)$, $(37,9,2)$ and $(56,11,2)$ are completely classified.
The largest set of parameters for which a biplane is known to exist is $(79,13,2)$.
In \cite{Ash}, Aschbacher constructed two biplanes with parameters $(79,13,2)$, dual to each other, with the full automorphism group of order $110$. 
So far, these are the only known biplanes with $79$ points. According to \cite[Theorem 3]{biplanes}, the order of the full automorphism group of any other biplane with parameters $(79,13,2)$ 
must be equal to $1$ or $3$. The next possible parameter triple for an existence of a symmetric $2$-$(v,k,2)$ design is $(121,16,2)$. 
There is no known example of a biplane with $121$ points. The following results are given in \cite{biplanes}.

\begin{result}(\cite[Lemma 14]{biplanes})\label{res_lm}
Let ${\mathcal D} =( {\mathcal P},{\mathcal B},{\mathcal I})$ be a biplane with parameters $(121,16,2)$. 
Then a Sylow $p$-subgroup of $Aut({\mathcal D})$ with $p \in \{5,7,11,13 \}$ is cyclic of order at most $p$. 
Moreover, the cyclic type of a $p$-element $x \in Aut({\mathcal D})$ on ${\mathcal P}$ and on a fixed block (if such exists) is given in Table \ref{res_tab}.
\end{result}

\begin{table}[h]  
\begin{center}
\begin{tabular}{c|c|c}
p & Cycle type on ${\mathcal P}$  & Cycle type on a fixed block $B$\\
\hline 
5& $1^1 5^{24}$ & $1^1 5^3$ \\
7& $1^2 7^{17}$ & $1^2 7^2$ \\
11& $11^{11}$ & $x$ fixes no blocks \\
13& $1^4 13^{9}$ & $1^3 13^1$ \\
\end{tabular}
\caption{\footnotesize Cycle type of $p$-elements $x$ with $p \in \{5,7,11,13 \}$ } \label{res_tab}
\end{center}
\end{table}

\begin{result} (\cite[Theorem 1, (b)]{biplanes}) \label{res_th}
Let ${\mathcal D}$ be a biplane with parameters $(121,16,2)$. Then $|Aut({\mathcal D})|$ divides $2^7\cdot3^2\cdot5\cdot7\cdot11\cdot13.$
\end{result}

The main goal of this paper is to give further restrictions on the order of an automorphism group of a biplane with parameters $(121,16,2)$. The main theorem of this paper is given below.

\begin{thm}[Main theorem] \label{main_th}
Let ${\mathcal D}$ be a biplane with parameters $(121,16,2)$. Then $|Aut({\mathcal D})|$ divides $2^7\cdot3^2\cdot5\cdot7$. Moreover, if $G\leq Aut({\mathcal D})$, then
\begin{enumerate}
 \item $|G| \notin \{14,35 \}$,
 \item if $|G|=21$, then $G\cong Frob_{21}$,
 \item  If  $|G|=10$, then $G\cong D_{10}$.
\end{enumerate}
\end{thm}

In order to prove Theorem \ref{main_th}, we will first show that it is not possible to construct an orbit matrix for an action of an automorphism of order $11$ or $13$ acting on a biplane with 
parameters $(121,16,2)$.
Then we study a possible action of an automorphism of order five or seven, and some small groups of order divisible by five or seven, on a biplane with parameters $(121,16,2)$.

For the construction of (partial) orbit matrices we used computer programs developed in scope of the PhD thesis of the second author (see \cite{doris}).

\section{Orbit matrices and automorphisms of order 11 or 13} \label{OM}

Let $ {\cal D}=( {\cal P},{\cal B},I) $ be a symmetric $ (v,k,\lambda ) $ design and $ G\leq {\rm Aut}({\cal D})$. The group action of $ G $ produces the same number
of point and block orbits (see \cite[Theorem 3.3]{land}). We denote that number by $ t $, the point orbits by ${\cal P}_{1},\ldots ,{\cal P}_{t} $,
the block orbits by $ {\cal B}_{1},\ldots ,{\cal B}_{t} $, and put $ |{\cal P}_{r}|=\omega _{r} $ and $ |{\cal B}_{i}|=\Omega _{i} $.
We denote by $ \gamma _{ir} $ the number of points of $ {\cal P}_{r} $ which are incident with a representative of the block orbit ${\cal B}_{i}$.
For these numbers the following equalities hold (see \cite{cep, crn-ruk, des-6, j}):
\begin{eqnarray}
\sum _{r=1}^{t}{\gamma }_{ir} & = & k \label{eq-1}\, ,\\ \label{def1}
\sum _{r=1}^{t}\frac{{\Omega }_{j}}{{\omega }_{r}}{\gamma }_{ir}{\gamma }_{jr} & = & \lambda {\Omega }_{j}+{\delta }_{ij}\cdot (k-\lambda )\label{eq-2}\, , \label{def2}
\end{eqnarray}
where ${\delta }_{ij}$ is Kronecker delta.  

\begin{defi} \label{def-orb-mat}
A $(t \times t)$-matrix $({\gamma}_{ir})$ with entries satisfying conditions $(\ref{eq-1})$ and $(\ref{eq-2})$ is called a $($block$)$ orbit matrix for
the parameters $(v,k, \lambda)$ and orbit lengths distributions $(\omega_{1}, \ldots ,\omega_{t})$, $(\Omega_{1}, \ldots ,\Omega_{t})$. 
\end{defi}

Hence, the $ith$ row of an orbit matrix satisfy the condition ($\ref{eq-1}$) and the condtion

\begin{eqnarray}
\sum _{r=1}^{t}\frac{{\Omega }_{i}}{{\omega }_{r}}{\gamma }_{ir}^2 & = & \lambda ({\Omega }_{i}-1)+k \label{eq-2a}\,.
\end{eqnarray}

Let us notice that in the case when an automorphism of prime order acts on a symmetric design, after a possible reordering of orbits, we have 
$(\omega_{1}, \ldots ,\omega_{t})=(\Omega_{1}, \ldots ,\Omega_{t})$.

\begin{defi}
An $(s \times t)$-matrix $({\gamma}_{ir})$, $s<t$, with entries satisfying condition $(\ref{eq-1})$, and the condition
\begin{eqnarray}
\sum _{r=1}^{t}\frac{{\Omega }_{j}}{{\omega }_{r}}{\gamma }_{ir}{\gamma }_{jr}  \leq \lambda {\Omega }_{j}+{\delta }_{ij}\cdot (k-\lambda )\label{eq-2'}\,. \label{pom}
\end{eqnarray}
is called a partial $($block$)$ orbit matrix for the parameters $(v,k, \lambda)$ and orbit lengths distributions $(\omega_{1}, \ldots ,\omega_{t})$, $(\Omega_{1}, \ldots ,\Omega_{t})$.  
\end{defi}

More information on orbit matrices of block designs and a construction of 2-designs using orbit matrices can be found in \cite{crn-ruk}.

\subsection{Orbit matrices for an automorphism of order 13} \label{OM_13}

According to Result \ref{res_lm}, the only possible orbit lengths distribution for an action of an automorphism of order $13$ on the points (and blocks) of a biplane with parameters $(121,16,2)$ is 
$(1,1,1,1,13,13,13,13,13,13,13,13,13)\sim 1^4 13^9$, and the cycle type on a fixed block is given by $1^3 13^1$. In other words, every fixed block is of type $(1,1,1,0,13,0,0,0,0,0,0,0,0)$, 
where the first four numbers represent the number of points in the orbits of length one incident with a fixed block, and the rest of them represent distribution of points in orbits of length $13$.  
One can easily check that the equalities ($\ref{eq-1}$) and ($\ref{eq-2a}$) hold for the type $(1,1,1,0,13,0,0,0,0,0,0,0,0)$ and $\Omega_i = 1$.

\begin{lem} \label{lema_aut13}
 Let ${\mathcal D}$ be a biplane with parameters $(121,16,2)$ and let $\rho$ be an automorphism of ${\mathcal D}$.  Then the order of $\rho$ is not equal to $13$.
\end{lem}
{\bf Proof}
The only possible action of $\rho$ is with four fixed points and nine orbits of length $13$ (Result \ref{res_lm}).
The corresponding orbit matrix $OM$ is a $13 \times 13$ matrix $({\gamma}_{ir})$ with entries satisfying conditions ($\ref{eq-1}$) and ($\ref{eq-2}$).
Let us determine the possible types of rows of $OM$. 

As stated above, $(1,1,1,0,13,0,0,0,0,0,0,0,0)$ is the unique type for a row corresponding to a fixed block. 
Taking into account the condition ($\ref{eq-2}$), without loss of generality we can assume that the first four rows (corresponding to the fixed blocks) of OM are given in Table \ref{OM13-fix}.

\begin{table}[h]
\begin{center} \begin{scriptsize}
\begin{tabular}{c|ccccccccccccc}
$OM$ & 1 & 1 & 1 & 1 & 13 & 13 & 13 & 13 & 13 & 13 & 13 & 13 & 13\\
\hline 
1 & 1 &  1 & 1 & 0 & 13 & 0 & 0& 0 & 0 & 0& 0& 0& 0\\
1 & 1 &  1 & 0 & 1 & 0 & 13 & 0& 0 & 0 & 0& 0& 0& 0\\
1 & 1 &  0 & 1 & 1 & 0 & 0 & 13& 0 & 0 & 0& 0& 0& 0\\
1 & 0 &  1 & 1 & 1 & 0 & 0 & 0& 13 & 0 & 0& 0& 0& 0\\
\end{tabular} \end{scriptsize}
\caption{\footnotesize Up to isomorphism, the unique partial orbit matrix corresponding to the four fixed blocks for an automorphism group of order $13$ acting on a biplane $(121,16,2)$} \label{OM13-fix}
\end{center}
\end{table}

We'll now determine possible types for rows corresponding to the block orbits of length $13$.
Since every fixed point should appear in exactly $16$ blocks, and all pairs of fixed points have already appeared twice in a common block, we need row types starting with $1,0,0,0$.
Each block belonging to an orbit of length 13 have to intersect each of the fixed blocks in exactly two points. So, we need a type having $1,0,0,0,2,1,1,1$ as its components. 
Because of the conditions ($\ref{eq-1}$) and ($\ref{eq-2}$), the rest of the positions should be occupied by five non-negative integers whose sum is equal to $10$ and the sum of their squares is $20$.
There is a unique solution satisfying those conditions, namely the type $(1,0,0,0,2,2,2,2,2,2,1,1,1)$. Rows of that type appear exactly four times in $OM$, and there is the unique partial orbit matrix corresponding to the blocks that contain a fixed point, which is given in Table \ref{OM13-fp}.   
\begin{table}[h]
\begin{center} \begin{scriptsize}
\begin{tabular}{c|ccccccccccccc}
$OM$ & 1 & 1 & 1 & 1 & 13 & 13 & 13 & 13 & 13 & 13 & 13 & 13 & 13\\
\hline 
1 & 1 &  1 & 1 & 0 & 13 & 0 & 0& 0 & 0 & 0& 0& 0& 0\\
1 & 1 &  1 & 0 & 1 & 0 & 13 & 0& 0 & 0 & 0& 0& 0& 0\\
1 & 1 &  0 & 1 & 1 & 0 & 0 & 13& 0 & 0 & 0& 0& 0& 0\\
1 & 0 &  1 & 1 & 1 & 0 & 0 & 0& 13 & 0 & 0& 0& 0& 0\\
\hline
13& 1& 0& 0& 0& 1& 1& 1& 2& 2& 2& 2& 2& 2\\
13& 0& 1& 0& 0& 1& 1& 2& 1& 2& 2& 2& 2& 2\\
13& 0& 0& 1& 0& 1& 2& 1& 1& 2& 2& 2& 2&2 \\
13& 0& 0& 0& 1& 2& 1& 1& 1& 2& 2& 2& 2& 2\\
\end{tabular} \end{scriptsize}
\caption{\footnotesize Up to isomorphism, the unique partial orbit matrix corresponding to the blocks that contain a fixed point for an automorphism group of order $13$ of a biplane $(121,16,2)$} 
\label{OM13-fp}
\end{center}
\end{table}

The ninth row of $OM$ corresponds to the orbit of length $13$. Since every fixed point already appears $16$ times, blocks from that orbit do not contain fixed points and intersect each fixed block in two points from the same orbit of length $13$ (the same holds for rows 10-13).
That gives the type having the numbers $0,0,0,0,2,2,2,2$ as its components. The rest of the positions should be occupied by non-negative integers whose sum is equal to $8$, satisfying that the sum of their squares is equal to $24$. The only possible type satisfying these conditions is $(0,0,0,0,4,2,2,2,2,2,2,0,0)$, and the unique partial orbit matrix corresponding to the first nine block orbits is given in Table \ref{OM13-ninth}.

\begin{table}[h]
\begin{center} \begin{scriptsize}
\begin{tabular}{c|ccccccccccccc}
$OM$ & 1 & 1 & 1 & 1 & 13 & 13 & 13 & 13 & 13 & 13 & 13 & 13 & 13\\
\hline 
1 & 1 &  1 & 1 & 0 & 13 & 0 & 0& 0 & 0 & 0& 0& 0& 0\\
1 & 1 &  1 & 0 & 1 & 0 & 13 & 0& 0 & 0 & 0& 0& 0& 0\\
1 & 1 &  0 & 1 & 1 & 0 & 0 & 13& 0 & 0 & 0& 0& 0& 0\\
1 & 0 &  1 & 1 & 1 & 0 & 0 & 0& 13 & 0 & 0& 0& 0& 0\\
\hline
13& 1& 0& 0& 0& 1& 1& 1& 2& 2& 2& 2& 2& 2\\
13& 0& 1& 0& 0& 1& 1& 2& 1& 2& 2& 2& 2& 2\\
13& 0& 0& 1& 0& 1& 2& 1& 1& 2& 2& 2& 2&2 \\
13& 0& 0& 0& 1& 2& 1& 1& 1& 2& 2& 2& 2& 2\\
\hline
13& 0& 0& 0& 0& 2& 2& 2& 2& 4& 2& 2& 0&0 \\
\end{tabular} \end{scriptsize}
\caption{\footnotesize Up to isomorphism, the unique partial orbit matrix corresponding to the first nine block orbits of an automorphism group of order $13$ acting on a biplane with parameters 
$(121,16,2)$} \label{OM13-ninth}
\end{center}
\end{table}

The $i$-th and the $j$-th row of type $(0,0,0,0,4,2,2,2,2,2,2,0,0)$ must satisfy the condition $\sum _{r=1}^{13}{\gamma }_{ir}{\gamma }_{jr}  =26,$ which is not possible because their product is always a multiple of $4$.
So, we can complete at most nine rows of an orbit matrix for an action of an automorphism of order $13$ on a biplane with parameters $(121,16,2)$. 
Therefore, a biplane with parameters $(121,16,2)$ cannot have an automorphism of order $13$.
{$\Box$}

\subsection{Orbit matrices for an automorphism of order 11} \label{OM_11}

According to Result \ref{res_lm}, an automorphism of order $11$ acts on a biplane with parameters $(121,16,2)$ without fixed points. So, if there is a biplane with parameters $(121,16,2)$ having 
an automorphism group $G$ of order 11, then $G$ acts with the orbit lengths distribution $(11,11,11,11,11,11,11,11,11,11,11)\sim 11^{11}$. 
In the following lemma we show that such an action is not possible.

\begin{lem} \label{lema_aut11}
 Let ${\mathcal D}$ be a biplane with parameters $(121,16,2)$ and let $\rho$ be an automorphism of ${\mathcal D}$.  Then the order of $\rho$ is not equal to $11$.
\end{lem}
{\bf Proof} Firstly, we will construct the possible row types for an orbit matrix of an automorphism of order $11$ acting on a biplane with parameters $(121,16,2)$ without fixed points and blocks. 
From the conditions ($\ref{eq-1}$) and ($\ref{eq-2a}$) it follows that $\sum _{r=1}^{11}{\gamma }_{ir}  =  16 $ and $\sum _{r=1}^{11} {\gamma }_{ir}^2  = 36 .$
It follows that $0 \leq {\gamma }_{ir} \leq 4 $, where number $4$ could appear at most once in a type.
Further calculations yield four possible types, which are given in Table \ref{OM11_types}.

\begin{table}[h]
\begin{center}
\begin{tabular}{c|c}
type &  \\
\hline 
1 & $(4,3,2,1,1,1,1,1,1,1,0)$\\
2 & $(4,2,2,2,2,1,1,1,1,0,0)$\\
3 & $(3,3,3,2,1,1,1,1,1,0,0)$\\
4 & $(3,3,2,2,2,2,1,1,0,0,0)$\\
\end{tabular}
\caption{\footnotesize The possible row types of an orbit matrix for an automorphism group of order $11$ acting on a biplane with parameters $(121,16,2)$} \label{OM11_types}
\end{center}
\end{table}

The corresponding orbit matrix $OM$ is a $11 \times 11$ matrix $({\gamma}_{ir})$ with entries satisfying the equation $\sum _{r=1}^{11}{\gamma }_{ir}{\gamma }_{jr}  =  22$, when $i\neq j$. 
There is to much possibilities to proceed without computers. With the help of a computer program developed in \cite{doris}, we get that
it is not possible to construct more than seven rows of an orbit matrix for an automorphism of order $11$ acting on a biplane with parameters $(121,16,2)$ without fixed points and blocks. 
Hence, a group of order 11 cannot act as an automorphism group of a biplane $(121,16,2)$.
{$\Box$}

 \section{Automorphisms of order 5 or 7} \label{aut}
 
Orbit lengths distributions and cycle types of fixed blocks for an automorphism of order $5$ or $7$ acting on a biplane with parameters $(121,16,2)$ are given by Result \ref{res_lm}.  
In this section we determine the possible row types of orbit matrices for automorphism groups of order 5 or 7 and construct partial orbit matrices corresponding to the blocks incident with at least one fixed point. 

\subsection{An action of an automorphism of order $7$}\label{aut7}

According to Result \ref{res_lm}, an automorphism of order seven acts on a biplane with parameters $(121,16,2)$ with two fixed points and $17$ orbits of length seven, and the unique type of a fixed block 
$B$ is given by  $(1,1,7,7,0,0,0,0,0,0,0,0,0,0,0,0,0) \sim 1^2 7^2$. 

Using similar arguments as in the proof of Lemma \ref{lema_aut13}, we established that up to isomorphism there are eight partial orbit matrices corresponding to orbits of an
automorphism of order 7 containing blocks incident with at least one fixed point of a biplane (121,16,2). These partial orbit matrices are given in Table \ref{OM7_fp}.  

\begin{table}[h]
\begin{center} \begin{scriptsize}
\begin{tabular}{c|c|ccccccccccccccccccc||c}
$OM_{1-8}$& & 1 & 1 & 7 & 7 & 7 & 7 & 7 & 7 & 7 & 7 & 7 & 7 & 7 &7 &7 &7 &7 &7 &7& row\\
\hline 
&1 & 1 &  1 & 7 & 7 & 0 & 0& 0 & 0 & 0& 0& 0& 0& 0& 0& 0& 0& 0& 0& 0&1\\
&1 & 1 &  1 & 0 & 0 & 7 & 7& 0 & 0 & 0& 0& 0& 0& 0& 0& 0& 0& 0& 0& 0&2\\
&7 & 1 &  0 & 1 & 0 & 1 & 0& 2 & 2 & 2& 1& 1& 1& 1& 1& 1& 1& 0& 0& 0&3\\
&7 & 1 &  0 & 0 & 1 & 0 & 1& 0 & 0 & 0& 1& 1& 1& 1& 1& 1& 1& 2& 2& 2&4\\
\hline \hline
$OM_1$&7 & 0 &  1 & 1 & 0 & 1 & 0& 2 & 1 & 1& 1& 1& 1& 1& 0& 0& 0& 2& 2& 1&5\\
&7 & 0 &  1 & 0 & 1 & 0 & 1& 1 & 1 & 0& 1& 1& 1& 1& 2& 2& 2& 0& 0& 1&6\\
\hline \hline
$OM_2$&7 & 0 &  1 & 1 & 0 & 1 & 0& 2 & 1 & 0& 2& 1& 1& 1& 1& 0& 0& 2& 1& 1&5\\
&7 & 0 &  1 & 0 & 1 & 0 & 1& 0 & 1 & 2& 0& 1& 1& 1& 1& 2& 2& 0& 1& 1&6\\
\hline \hline
$OM_3$&7 & 0 &  1 & 1 & 0 & 1 & 0& 2 & 1 & 0& 1& 1& 1& 1& 1& 1& 0& 2& 2& 0&5\\
&7 & 0 &  1 & 0 & 1 & 0 & 1& 0 & 1 & 2& 1& 1& 1& 1& 1& 1& 2& 0& 0& 2&6\\
\hline \hline
$OM_4$&7 & 0 &  1 & 1 & 0 & 1 & 0& 1 & 1 & 1& 2& 2& 1& 1& 0& 0& 0& 2& 1& 1&5\\
&7 & 0 &  1 & 0 & 1 & 0 & 1& 1 & 1 & 1& 0& 0& 1& 1& 2& 2& 2& 0& 1& 1&6\\
\hline \hline
$OM_5$&7 & 0 &  1 & 1 & 0 & 0 & 1& 2 & 2 & 0& 1& 1& 1& 1& 1& 0& 0& 2& 1& 1&5\\
&7 & 0 &  1 & 0 & 1 & 1 & 0& 0 & 0 & 2& 1& 1& 1& 1& 1& 2& 2& 0& 1& 1&6\\
\hline \hline
$OM_6$&7 & 0 &  1 & 1 & 0 & 0 & 1& 2 & 1 & 1& 2& 1& 1& 1& 0& 0& 0& 2& 1& 1&5\\
&7 & 0 &  1 & 0 & 1 & 1 & 0& 0 & 1 & 1& 0& 1& 1& 1& 2& 2& 2& 0& 1& 1&6\\
\hline \hline
$OM_7$&7 & 0 &  1 & 1 & 0 & 0 & 1& 2 & 1 & 0& 2& 1& 1& 1& 1& 1& 0& 2& 1& 0&5\\
&7 & 0 &  1 & 0 & 1 & 1 & 0& 0 & 1 & 2& 0& 1& 1& 1& 1& 1& 2& 0& 1& 2&6\\
\hline \hline
$OM_8$&7 & 0 &  1 & 1 & 0 & 0 & 1& 1 & 1 & 1& 2& 2& 2& 1& 0& 0& 0& 1& 1& 1&5\\
&7 & 0 &  1 & 0 & 1 & 1 & 0& 1 & 1 & 1& 0& 0& 0& 1& 2& 2& 2& 1& 1& 1&6\\
\end{tabular} \end{scriptsize}
\caption{\footnotesize Partial orbit matrices corresponding to the blocks containing fixed points for an automorphism group of order seven acting on a biplane with parameters $(121,16,2)$} \label{OM7_fp}
\end{center}
\end{table}

Further, using the equalities (\ref{eq-1}) and (\ref{eq-2a}), we established that the four types given in Table \ref{OM7_types} are the only row types for orbits containing blocks that are not incident with any of the fixed points.

\begin{table}[h]
\begin{center}
\begin{tabular}{c|c}
type &  \\
\hline 
1 & $(0,0,4,1,1,1,1,1,1,1,1,1,1,1,1,0,0,0,0)$\\
2 & $(0,0,3,3,1,1,1,1,1,1,1,1,1,1,0,0,0,0,0)$\\
3 & $(0,0,3,2,2,2,1,1,1,1,1,1,1,0,0,0,0,0,0)$\\
4 & $(0,0,2,2,2,2,2,2,1,1,1,1,0,0,0,0,0,0,0)$\\
\end{tabular} 
\caption{\footnotesize The possible types for a row of an orbit matrix for an automorphism of order seven acting on a biplane with parameters $(121,16,2)$ corresponding to blocks not incident with any fixed point} \label{OM7_types}
\end{center}
\end{table}

\begin{rmk} 
Because of the large number of possibilities, it is necessary to involve computers in a construction of further orbits.
In that way we obtained $2999$ partial orbit matrices by adding just one row to partial orbit matrices from Table \ref{OM7_fp}, and more than one million possibilities for partial orbit matrices consisting of eight rows. The construction of full orbit matrices is out of our reach.
\end{rmk}

\subsection{An action of an automorphism of order $5$} \label{aut5}

An automorphism of order five acts on a biplane with parameters $(121,16,2)$ with one fixed point and $24$ orbits of length five, and the unique type of a fixed block $B$ is given by $1^1 5^3$. 
The unique partial orbit matrix corresponding to the blocks that contain fixed point, obtained in a similar way as described in the proof of Lemma \ref{lema_aut13}, is given in the Table \ref{OM5_fp}.

\begin{table}[h]
\begin{center} \begin{tiny}
\begin{tabular}{c|ccccccccccccccccccccccccc}
$OM$ & 1 & 5 & 5 & 5 & 5 & 5 & 5 & 5 & 5 & 5 & 5 & 5 & 5 &5 &5 &5 &5 &5 &5& 5& 5& 5& 5& 5& 5\\
\hline 
1 & 1 &  5 & 5 & 5 & 0 & 0& 0 & 0 & 0& 0& 0& 0& 0& 0& 0& 0& 0& 0& 0&0 &0 &0 &0 &0 &0\\
5 & 1 &  1 & 0 & 0 & 2 & 2& 1 & 1 & 1& 1& 1& 1& 1& 1& 1& 1& 0& 0& 0&0 &0 &0 &0 &0 &0\\
5 & 1 &  0 & 1 & 0 & 0 & 0& 1 & 1 & 1& 1& 1& 0& 0& 0& 0& 0& 2& 2& 1&1 &1 &1 &1 &0 &0\\
5 & 1 &  0 & 0 & 1 & 0 & 0& 0 & 0 & 0& 0& 0& 1& 1& 1& 1& 1& 0& 0& 1&1 &1 &1 &1 &2 &2\\
 
\end{tabular} \end{tiny}
\caption{\footnotesize The partial orbit matrix corresponding to the blocks containing fixed point for an automorphism of order five acting on a biplane with parameters $(121,16,2)$} \label{OM5_fp}
\end{center}
\end{table}

Using the equalities (\ref{eq-1}) and (\ref{eq-2a}) we established that the two types given in Table \ref{OM5_types} are the only row types for orbits of an automorphism group of order five containing blocks that are not incident with any of the fixed points.  

\begin{table}[h]
\begin{center}
\begin{tabular}{c|c}
type &  \\
\hline 
1 & $(0,3,2,1,1,1,1,1,1,1,1,1,1,1,0,0,0,0,0,0,0,0,0,0,0)$\\
2 & $(0,2,2,2,2,1,1,1,1,1,1,1,1,0,0,0,0,0,0,0,0,0,0,0,0)$\\
\end{tabular} 
\caption{\footnotesize The possible types for a row of an orbit matrix for an automorphism of order five acting on a biplane with parameters $(121,16,2)$ corresponding to blocks not incident with the fixed point} \label{OM5_types}
\end{center}
\end{table}

\begin{rmk} It is necessary to involve computers in a further construction.
In that way we obtained $346$ partial orbit matrices by adding just one row to the partial orbit matrix from Table \ref{OM5_fp}, and $758662$ partial orbit matrices with six rows.
The construction of the orbit matrices is out of our reach.
\end{rmk}

\section {Small groups of composite order divisible by 5 or 7}

In this section we study a possible action of automorphism groups of order 10, 14, 15, 21 and 35 on a biplane with parameters (121,16,2).

Let a group $G$ acts on a nonempty set $\Omega$ and $H \unlhd G$. Then each $G$-orbit of $\Omega$ decomposes to $H$-orbits of the same size and the group $G/H$ acts transitively on these $H$-orbits 
(see \cite{crn-ruk}). We'll use this fact in the sequel. Note that automorphism groups of order 10, 14, 15, 21 and 35 are direct or semidirect products of cyclic groups. 

\begin{lem} \label{lema_aut35}
Let ${\mathcal D}$ be a biplane with parameters $(121,16,2)$. If $G\leq Aut({\mathcal D})$, then $|G| \neq 35$.
\end{lem}
{\bf Proof}
The only group of order $35$ is the cyclic group $Z_{35} \cong Z_5 \times Z_7$. According to Result \ref{res_lm}, an automorphism of order five fixes one point, and an automorphism of order seven fixes two points. This is not possible, because an automorphism of order seven fixes every point belonging to a $G$-orbit of length one or five (and an automorphism of order five fixes every point belonging to a 
$G$-orbit of length one or seven). 
{$\Box$}

\begin{lem} \label{lema_aut21}
Let ${\mathcal D}$ be a biplane with parameters $(121,16,2)$. If $G \leq Aut({\mathcal D})$ and $|G| = 21$, then $G \cong Frob_{21} \cong Z_7 : Z_3$.
\end{lem}
{\bf Proof}
There are two groups of order $21$. One of them is the Frobenius group $Frob_{21}$ which is isomorphic to the semidirect product $Z_7 : Z_3$, and the other is the cyclic group $Z_{21}\cong Z_3 \times  Z_7$.
An automorphism of order three fixes one or seven points (see \cite[Lemma 13]{biplanes}), and an automorphism of order seven fixes two points.

If $G \cong Z_{21}$, then an automorphism of order seven fixes every point belonging to a $G$-orbit of length one or three and an automorphism of order three fixes every point belonging to a 
$G$-orbit of length one or seven. Hence, $G$ cannot act in a way that an automorphism of order three fixes one or seven points, and an automorphism of order seven fixes two points.
{$\Box$}

If $G \cong Frob_{21}$, then an automorphism of order seven fixes every point belonging to a $G$-orbit of length three, but an automorphism of order three fixes just one point from a $G$-orbit of length
seven. There is a unique orbit lengths distribution for an action of the group $Frob_{21}$ on ${\mathcal D}$, namely the distribution $(1,1,7,7,7,7,7,21,21,21,21)$.
Further, there is the unique orbit matrix for that distribution, which is given in Table \ref{OM21}.

\begin{table}[h]
\begin{center}  
\begin{tabular}{c|ccccccccccc}
$OM$ & 1 & 1 & 7 & 7 & 7 & 7 & 7 & 21 & 21 & 21 & 21\\
\hline 
1 & 1 &  1 & 7 & 7 & 0 & 0 & 0& 0 & 0 & 0& 0\\
1 & 1 &  1 & 0 & 0 & 7 & 7 & 0& 0 & 0 & 0& 0\\
 
7& 1& 0& 1& 0& 1& 0& 1& 6& 3& 3& 0\\
7& 1& 0& 0& 1& 0& 1& 1& 0& 3& 3& 6\\
7& 0& 1& 1& 0& 0& 1& 1& 3& 6& 0& 3\\
7& 0& 1& 0& 1& 1& 0& 1& 3& 0& 6& 3\\
7& 0& 0& 1& 1& 1& 1& 0& 6& 0& 0& 6\\

21& 0& 0& 2& 0& 1& 1& 2& 2& 2& 3& 3\\
21& 0& 0& 0& 2& 1& 1& 2& 3& 3& 2& 2\\
21& 0& 0& 1& 1& 2& 0& 0& 2& 4& 3& 3\\
21& 0& 0& 1& 1& 0& 2& 0& 3& 3& 4& 2\\
\end{tabular}  
\caption{\footnotesize The unique orbit matrix for an action of a group of order $21$ on a biplane $(121,16,2)$} \label{OM21}
\end{center}
\end{table}

Because of the large number of possibilities, construction of biplanes $(121,16,2)$ corresponding to the orbit matrix $OM$ is out of our range. 

\begin{rmk}
The group $Z_{15}$ is the only group of order 15. If $G\cong Z_{15}$ acts on a biplane $(121,16,2)$, then it acts with the orbit lengths distribution $(1,15,15,15,15,15,15,15,15)$. 
Up to isomorphism there are exactly six orbit matrices for such an action of $G$. These orbit matrices are given in Appendix.
\end{rmk}

\begin{lem} \label{lema_aut14}
Let ${\mathcal D}$ be a biplane with parameters $(121,16,2)$. If $G\leq Aut({\mathcal D})$, then $|G| \neq 14$.
\end{lem}
{\bf Proof}
There are two groups of order $14$, the dihedral group $D_{14} \cong Z_7 : Z_2$ and the cyclic group $Z_{14}\cong Z_2 \times  Z_7.$
According to \cite{z14}, a biplane with parameters $(121,16,2)$ with an automorphism group $Z_{14}$ does not exist.

Let $G$ be isomorphic to $D_{14}$. An involution fixes nine or thirteen points (see \cite[Lemma 10]{biplanes}), and the number of fixed points for an automorphism of order seven is two.
Further, an involution acts on $G$-orbit of length seven fixing one point.
Therefore, the possible orbit lengths distributions for $D_{14}$ acting on a biplane with parameters (121,16,2) are $(1,1,7,7,7,7,7,7,7,14,14,14,14,14)$ and $(1,1,7,7,7,7,7,7,7,7,7,7,7,14,14,14)$.
Using computer programs developed in \cite{doris}, we obtained that it is not possible to construct orbit matrices for these orbit lengths distributions.
{$\Box$}

\begin{lem} \label{lema_aut10}
Let ${\mathcal D}$ be a biplane with parameters $(121,16,2)$. If $G \leq Aut({\mathcal D})$ and $|G| = 10$, then $G$ is isomorphic to the dihedral group $D_{10}$.
\end{lem}
{\bf Proof}
There are two groups of order $10$, the dihedral group $D_{10} \cong Z_5 : Z_2$ and the cyclic group $Z_{10}\cong Z_2 \times  Z_5$.
An involution acting on ${\mathcal D}$ fixes $9$ or $13$ points from ${\mathcal P}$ and an automorphism of order five acting on ${\mathcal D}$ fixes one point.

Let $G$ be isomorphic to $Z_{10}$. Then every involution in $G$ fixes every point from $G$-orbit of length one or five, and an automorphism of order five fixes every point from $G$-orbit of length one or two. Such an action on the point set of ${\mathcal D}$ is not possible.

Let $G$ be isomorphic to the dihedral group $D_{10}$. An involution fixes one point from $G$-orbit of length five, and an automorphism of order five fixes every point from $G$-orbit of length one or two.
The only possible orbit lengths distributions are
\begin{enumerate} 
\item $(1,5,5,5,5,5,5,5,5,10,10,10,10,10,10,10,10)$, 
\item $(1,5,5,5,5,5,5,5,5,5,5,5,5,10,10,10,10,10,10)$.
\end{enumerate}
For the first distribution we obtained $24$ orbit matrices. In the case of the second orbit lengths distribution there are to many possibilities and we were not able to build the corresponding orbit matrices, or to show that the orbit matrices does not exist.
{$\Box$}

\bigskip

The statement of Main theorem follows from Lemmas \ref{lema_aut13}, \ref{lema_aut11}, \ref{lema_aut35}, \ref{lema_aut21}, \ref{lema_aut14} and \ref{lema_aut10}.

\vspace*{0.5cm}
\begin{center}{\bf Acknowledgement}\end{center}
This work has been fully supported by {\rm C}roatian Science Foundation under the project 6732.

\newpage

\section*{Appendix} \label{Ap_A}

Below we give the orbit matrices for a group $G\cong Z_{15}$.

\bigskip

 \begin{footnotesize} \noindent  
 
\noindent 
$ \begin{array}{c|ccccccccc}
OM_1 & 1& 15 & 15 & 15 & 15 & 15 & 15 & 15 & 15 \\
\hline 
1 & 1 &  15 & 0 & 0 & 0 & 0 & 0& 0 & 0\\
 
15& 1& 1& 2& 2& 2& 2& 2& 2& 2\\

15& 0& 2& 5& 2& 2& 2& 1& 1& 1\\

15& 0& 2& 2& 3& 1& 0& 4& 3& 1\\
15& 0& 2& 2& 1& 0& 3& 1& 4& 3\\
15& 0& 2& 2& 0& 3& 1& 3& 1& 4\\

15& 0& 2& 1& 4& 3& 1& 0& 2& 3 \\
15& 0& 2& 1& 3& 1& 4& 3& 0& 2 \\
15& 0& 2& 1& 4& 3& 1& 2& 3& 0 \\
 \end{array} $   \hfill \quad
$ \begin{array}{c|ccccccccc}
OM_1 & 1& 15 & 15 & 15 & 15 & 15 & 15 & 15 & 15 \\
\hline 
1 & 1 &  15 & 0 & 0 & 0 & 0 & 0& 0 & 0\\
 
15& 1& 1& 2& 2& 2& 2& 2& 2& 2\\

15& 0& 2& 5& 2& 2& 2& 1& 1& 1\\

15& 0& 2& 2& 3& 1& 0& 4& 3& 1\\
15& 0& 2& 2& 1& 0& 3& 1& 4& 3\\
15& 0& 2& 2& 0& 3& 1& 3& 1& 4\\

15& 0& 2& 1& 4& 1& 3& 0& 3& 2 \\
15& 0& 2& 1& 3& 4& 1& 2& 0& 3 \\
15& 0& 2& 1& 1& 3& 4& 3& 2& 0 \\
\end{array} $  
 
\bigskip

\noindent 
$ \begin{array}{c|ccccccccc}
OM_3 & 1& 15 & 15 & 15 & 15 & 15 & 15 & 15 & 15 \\
\hline 
1 & 1 &  15 & 0 & 0 & 0 & 0 & 0& 0 & 0\\
 
15& 1& 1& 2& 2& 2& 2& 2& 2& 2\\

15& 0& 2& 5& 2& 2& 2& 1& 1& 1\\

15& 0& 2& 2& 1& 0& 3& 1& 4& 3\\
15& 0& 2& 2& 3& 1& 0& 3& 1& 4\\
15& 0& 2& 2& 0& 3& 1& 4& 3& 1\\

15& 0& 2& 1& 4& 3& 1& 0& 2& 3 \\
15& 0& 2& 1& 1& 4& 3& 3& 0& 2 \\
15& 0& 2& 1& 3& 1& 4& 2& 3& 0 \\
\end{array} $  \hfill \quad
$ \begin{array}{c|ccccccccc}
OM_4 & 1& 15 & 15 & 15 & 15 & 15 & 15 & 15 & 15 \\
\hline 
1 & 1 &  15 & 0 & 0 & 0 & 0 & 0& 0 & 0\\
 
15& 1& 1& 2& 2& 2& 2& 2& 2& 2\\

15& 0& 2& 5& 2& 2& 2& 1& 1& 1\\

15& 0& 2& 2& 3& 1& 0& 4& 3& 1\\
15& 0& 2& 2& 0& 3& 1& 1& 4& 3\\
15& 0& 2& 2& 1& 0& 3& 3& 1& 4\\

15& 0& 2& 1& 4& 1& 3& 2& 0& 3 \\
15& 0& 2& 1& 3& 4& 1& 3& 2& 0 \\
15& 0& 2& 1& 1& 3& 4& 0& 3& 2 \\
\end{array} $  
 
\bigskip

\noindent 
$ \begin{array}{c|ccccccccc}
OM_5 & 1& 15 & 15 & 15 & 15 & 15 & 15 & 15 & 15 \\
\hline 
1 & 1 &  15 & 0 & 0 & 0 & 0 & 0& 0 & 0\\
 
15& 1& 1& 2& 2& 2& 2& 2& 2& 2\\

15& 0& 2& 4& 4& 1& 1& 1& 1& 1\\

15& 0& 2& 3& 1& 0& 4& 3& 2& 1\\
15& 0& 2& 3& 0& 3& 2& 1& 1& 4\\
15& 0& 2& 2& 1& 3& 0& 4& 3& 1\\

15& 0& 2& 1& 3& 0& 1& 2& 3& 4 \\
15& 0& 2& 1& 2& 3& 3& 0& 4& 1 \\
15& 0& 2& 0& 3& 3& 3& 3& 0& 2 \\
\end{array} $  \hfill \quad
$ \begin{array}{c|ccccccccc}
OM_6 & 1& 15 & 15 & 15 & 15 & 15 & 15 & 15 & 15 \\
\hline 
1 & 1 &  15 & 0 & 0 & 0 & 0 & 0& 0 & 0\\
 
15& 1& 1& 2& 2& 2& 2& 2& 2& 2\\

15& 0& 2& 4& 4& 1& 1& 1& 1& 1\\

15& 0& 2& 3& 1& 0& 4& 3& 2& 1\\
15& 0& 2& 3& 0& 3& 1& 1& 4& 2\\
15& 0& 2& 2& 1& 3& 1& 4& 0& 3\\

15& 0& 2& 1& 3& 0& 1& 2& 3& 4 \\
15& 0& 2& 1& 2& 3& 4& 0& 1& 3 \\
15& 0& 2& 0& 3& 3& 2& 3& 3& 0 \\
\end{array} $
\end{footnotesize}

\end{document}